\documentclass[12pt]{article}
\usepackage[pctex32]{graphics}
\usepackage{amssymb}
\usepackage{latexsym}
\usepackage{epsfig} 
\usepackage{multirow} 
\usepackage{pstcol}
\usepackage{float}
\usepackage{epstopdf}
\usepackage{mathrsfs}
\usepackage{soul}

\usepackage{amsmath}
\usepackage{amssymb}
\usepackage{relsize}
\usepackage{multirow}

\usepackage{graphicx}
\usepackage{multicol} 

\usepackage{times}
\usepackage{epstopdf}

\newtheorem{definition}{Definition}[section]
\newtheorem{theorem}[definition]{Theorem}

\definecolor{Red}{rgb}{1,0.,0.}
\usepackage{fancyhdr}
\usepackage{setspace}
\usepackage{hyperref}
\usepackage{mathtools}
\usepackage{algorithm}
\usepackage[dvipsnames]{xcolor}
\usepackage{algorithmic}
\usepackage{placeins}
\newcommand{\R}{{\mathbb R}}

\newcommand{\bbE}{{\mathbb E}}

\newcommand{\mP}{{\mathsf P}}
\newcommand{\mT}{{\mathsf T}}
\newcommand{\mA}{{\mathsf A}}

\newcommand{\mI}{{\mathsf I}}
\newcommand{\mC}{{\mathsf C}}
\newcommand{\mD}{{\mathsf D}}
\newcommand{\mH}{{\mathsf H}}

\newcommand{\mV}{{\mathsf V}}
\newcommand{\mU}{{\mathsf U}}

\newcommand{\mG}{{\mathsf G}}

\newcommand{\mZ}{{\mathsf Z}}

\newcommand{\mOmega}{{\mathsf \Omega}}

\title{Spotlight inversion by orthogonal projections}
\author{D Calvetti$^1$ \and N Hyv\"{o}nen$^2$ \and V Kolehmainen$^3$ \and E Somersalo$^1$}
\date{$^1$ Case Western Reserve University, Cleveland, OH, USA\\
$^2$ Aalto University, Espoo, Finland\\
$^3$ Department of Technical Physics, University of Eastern Finland, Kuopio, Finland}

\begin{document}
\maketitle
\begin{abstract}
Many computational problems involve solving a linear system of equations, although only a subset of the entries of the solution are needed. In inverse problems, where the goal is to estimate unknown parameters from indirect noisy observations, it is not uncommon that the forward model linking the observed variables to the unknowns depends on variables that are not of primary interest, often referred to as nuisance parameters. In this article, we consider linear problems, and propose a novel projection technique to eliminate, or at least mitigate, the contribution of the nuisance parameters in the model.  We refer to this approach as spotlight inversion, as it allows to focus on only the portion of primary interest of the unknown parameter vector, leaving the uninteresting part in the shadow. The viability of the approach is illustrated with two computed examples, one where it works as model reduction for a finite element approximation of an elliptic PDE, the other amounting to local fanbeam X-ray tomography, spotlighting the region of interest that is part of the full target.
\end{abstract}

\section{Introduction}

    In numerous applications in science and engineering it is   necessary to solve systems of linear equations with interest in only portions of the solution vector. These problems  can be formulated as follows: Given a linear model
\[
 b = \mA x + \varepsilon,
\]
where $\mA\in\R^{m\times n}$ is a known matrix, $b$ is a given vector, and $\varepsilon$ is possibly vanishing random additive noise,   estimate {\em selected components} of the vector $x\in\R^n$. Without loss of generality, we can permute the entries of $x$ so that the quantity of interest consists of its first $n_1$ components, $n_1<n$. By partitioning the matrix $\mA$ accordingly, we can reformulate the problem as
\begin{equation}\label{original}
 b = \mA_1 x_1 + \mA_2 x_2 + \varepsilon,
\end{equation}
where $\mA_j \in \R^{m\times n_j}$, $j=1,2$, are known matrices, $n_1+n_2 = n$,
and $x_j\in\R^{n_j}$ are vectors,
\[
 \mA = \left[\begin{array}{cc} \mA_1 & \mA_2\end{array}\right], \quad x = \left[\begin{array}{c} x_1 \\ x_2\end{array}\right].
\]

Applications that are amenable to this type of formulation include inverse problems ($\varepsilon \neq 0$) in which $b$ represents a noisy observation, and forward problems ($\varepsilon = 0$) where the matrix $\mA$ is typically a square and invertible stiffness matrix arising from an elliptic partial differential equation and $b$ is a known load vector. In the latter case, one may be interested only in a selected subset of the components of the solution $x$,  discarding  the components of no relevance. In inverse problems, computational tasks of the form (\ref{original}) arise in  many different contexts, including the following important applications. 

{\em Local tomography:} In local X-ray tomography, the measured X-ray attenuation data depends on the absorption properties of the body extending over the entire field of view of the imaging system, although only a portion of the domain, often referred to as region of interest (ROI) corresponding to $x_1$, may be of interest. Local tomography commonly arises in X-ray tomography applications, either because of limited detector size, e.g., in dental Cone Beam Computed Tomography (CBCT) or because of the need to  reduce the X-ray dose outside of the ROI, e.g., in cardiac imaging.

{\em Radar imaging:} In radar imaging, the signal comprises the portion of interest coming from the main lobe of the radar, corresponding to $x_1$, and the clutter from the radar sidelobes that pick up the echo of trees, ground etc.,  represented by $x_2$.

{\em Image deblurring:} A blurred image typically accounts for the contribution from both the field of view and leak data from the outside region. Inappropriate modeling or ignoring of the latter may cause boundary artifacts that are often compensated for by means of artificial boundary conditions.

{\em MEG/EEG based brain imaging:} The search for the onset foci of epileptic seizures often can be restricted to a portion of the brain; Likewise, localization of active brain regions constituting a task-related, or task-negative connected network may be limited to selected brain regions, however the data account for contributions from the entire brain. 

{\em Functional MRI imaging:} Identification of a stimulus-related cerebral activity of interest typically may be circumscribed to a known part of the brain, and the 
dynamics of the rest of the brain and other tissues is of less or no interest.

In statistics, the portion of the unknown  $x_2$ in (\ref{original}) not of interest is referred to as {\em nuisance parameters} while the radar community refer to the contribution $\mA_2 x_2$ to the noiseless signal $\mA x$ as {\em clutter}. In this paper we adopt the terminology from both these two communities. 

In the inverse problems literature, the  separating the portion of interest from the clutter has been addressed in different ways. The tools for nuisance parameter elimination in classical statistics include considerations of sufficient and ancillary statistics and likelihood ratios \cite{basu2010elimination,linnik2008statistical}.
Bayesian statistics addresses the problem typically through integrated likelihood methods \cite{berger1999integrated} or marginalization over the nuisance parameter, often requiring MCMC sampling \cite{doucet2002marginal}. An example of optimization-based methods for estimating nuisance parameters in inverse problems can be found in
\cite{aravkin2012estimating}.
In a na\"{\i}ve approach, the nuisance parameter $x_2$ can be fixed to a plausible or representative value, or ignored by setting $x_2=0$, typically introducing a bias that in the Bayesian setting can be compensated by appropriate error modeling using the prior distribution of $x_2$ \cite{kaipio2007statistical}.

In this article, we put the spotlight of the inversion process on the portion $x_1$ of the unknown using a linear algebra-based approach exploiting suitable  orthogonal projections. Let 
\[
 \mP:\R^m\to  {\mathcal R}(\mA_2), \quad \mP^\perp:\R^m\to  {\mathcal R}(\mA_2)^\perp = {\mathcal N}(\mA_2^\mT),
\] 
denote the pair of orthogonal projectors onto the range ${\mathcal R}(\mA_2)\subset\R^m$ of the matrix $\mA_2$ and onto the the null space ${\mathcal N}(\mA_2^\mT)$ of the matrix $\mA_2^\mT$, respectively . 
Assuming that $m>n_2$, the subspace $ {\mathcal N}(\mA_2^\mT)$ is non-trivial. 
Applying $\mP^\perp$ to both sides of (\ref{original}) and observing that $\mP^\perp\mA_2 x_2 =0$ yields 
\begin{equation}\label{projected}
b' =  \mA_1' x_1 + \varepsilon',
\end{equation}
where $b' = \mP^\perp b$,  $\mA_1' = \mP^\perp\mA_1$, $\varepsilon' = \mP^\perp \varepsilon$. The projected problem, which does not contain the nuisance parameters can be solved as a least squares problem, equipped with a standard regularization technique if needed, or alternatively it can be analyzed in the framework of Bayesian inverse problems. 

In section~\ref{sec:forward}, we discuss the spotlighting for well-posed problems, starting with an $n\times n$ invertible matrix $\mA$, and we show that the approach reduces the original problem into an overdetermined least squares problem that recovers the selected components $x_1$ exactly.

In section~\ref{sec:bayesian}, the spotlighting  idea is applied to inverse problems where the dimensionality of the nuisance parameter does not exceed that of the data, $m>n_2$ In particular, we show that if the inverse problem is considered in a Bayesian setting, the proposed projection method for reducing the problem dimensionality does not correspond to the standard dimension reduction by marginalization in the Bayesian framework of inverse problems, but can be interpreted as a result of a limiting process as the underlying prior becomes uninformative. We emphasize that the clutter removal method by projection does not require or assume a Bayesian interpretation of the problem.

The assumption that the clutter is restricted to a proper subspace of $\R^m$, is automatically satisfied in the case of well-posed invertible problems $m=n$ as well as in overdetermined and formally determined {inverse} problems where  $m\geq n>n_2$. In inverse problems, the condition $m>n_2$ is not always guaranteed, and may not be satisfied for commonly encountered underdetermined problems. In such cases, by replacing the exact projector $\mP^\perp$ by a low-rank approximation,  spotlighting may still be a viable clutter reduction method, as will be discussed in \ref{sec:partial}.  Extensions to non-linear inverse problems are briefly discussed in section~\ref{sec:outlook}.

In the inverse problems literature, dimension reduction techniques have been discussed in the Bayesian context, see, e.g., \cite{spantini2015optimal,zahm2022certified}, based on the analysis of the informative directions of the likelihood with respect to the prior information. Unlike the cited work, which assume a Bayesian framework, our dimension reduction for inverse problems is purely linear algebraic and does not require a stochastic extension of the model. In particular, we need not assume any prior model for the unknown, and in the noiseless case $\varepsilon = 0$, even the likelihood is a futile concept.
A projection technique similar to the one proposed has been proposed in statistical regression analysis and in particular, in econometrics, referred to as partial regression, or (Yule) -Frisch-Waugh-Lovell ((Y)FWL) theorem 
\cite{davidson1993estimation} dating back to the early 20th century; see \cite{basu2023yule} for a historical account. The present article aims at bringing the idea into a computational linear algebraic context for both well-posed linear systems and ill-posed inverse problems.

\section{Spotlight projections for Model Reduction in Noiseless Case}\label{sec:forward}

Considering a noiseless model (\ref{original}), where $\varepsilon = 0$,
and $\mA\in\R^{n\times n}$ is an invertible matrix, partitioned as
\[
 \mA = \left[\begin{array}{cc} \mA_1 & \mA_2\end{array}\right], \quad \mA_1\in\R^{n\times k}, \quad \mA_2 \in\R^{n\times(n-k)},
\]
where $1\leq k<n$. 
Let $\mP: \R^n \to {\mathcal R}(\mA_2)$ and $\mP^\perp: \R^n \to {\mathcal R}(\mA_2)^\perp$ be the orthogonal projectors onto the range of $\mA_2$ and its orthogonal complement, respectively. Consider the equation
\begin{equation}\label{full}
 \mA x  = b,
\end{equation}
and its projected version onto the orthogonal complement of the range of $\mA_2$. 
\begin{equation}\label{projected}
 \mA_1' x_1 = b', \quad \mA_1' = \mP^\perp \mA_1, \quad b' = \mP^\perp b.
\end{equation}  
We start by stating and proving the following result.

 \begin{theorem}
 If $\mA\in \R^{n\times n}$ is invertible, and $\mA x = b$, then $x_1\in\R^k$ can be uniquely recovered by solving the projected least squares problem 
 \[
  \mA_1' x_1 = b'.
 \] 
 \end{theorem}

{\bf Proof.} First, observe that by the invertibility of $\mA$, its columns are linearly independent, hence
\[
 {\rm dim}\big({\mathcal R}(\mA_1)\big) = k, \quad {\rm dim}\big({\mathcal R}(\mA_2)\big) = n-k.
\]
Moreover,  
\[
 \R^n = {\mathcal R}(\mA_1) \oplus {\mathcal R}(\mA_2),
\]
that is, the space $\R^n$ is the direct, not necessarily orthogonal sum of the range of $\mA_1$ and the range of $\mA_2$. To prove this, we need to show that
\begin{enumerate}
\item $ \R^n = {\mathcal R}(\mA_1) + {\mathcal R}(\mA_2)$, i.e, every vector in $\R^n$ can be expressed as the sum of a vector in ${\mathcal R}(\mA_1) $ and a vector in ${\mathcal R}(\mA_2)$.
\item  $ {\mathcal R}(\mA_1) \cap {\mathcal R}(\mA_2) = \{0\}$, implying that the representation is unique.
\end{enumerate} 
To prove the first claim, observe that the columns of $\mA$ are linearly independent, hence a basis for $\R^n$, thus every $y\in\R^n$ admits a unique representation of the form 
\[
 y = \mA x = \left(\sum_{j=1}^k x_j a^{(j)}\right) +   \left(\sum_{j=k+1}^n x_j a^{(j)}\right),
\]
where $x = \mA^{-1} y$. 
To show the second claim, observe that if $z \in {\mathcal R}(\mA_1) \cap {\mathcal R}(\mA_2)$, it can be written in two ways,
\[
 z =   \left(\sum_{j=1}^k z_j a^{(j)}\right)  =    \left(\sum_{j=k+1}^n z_j a^{(j)}\right),
\]
which implies that the vector
\[
 \widetilde z = \left[\begin{array}{r} z_1 \\ \vdots \\ z_k \\ -z_{k+1} \\ \vdots \\ -z_n\end{array}\right],
\]
satisfies $\mA \widetilde z = 0$, so by the invertibility of $\mA$, $\widetilde z = 0$ and hence $z=0$.  

 To show that the matrix $\mA_1'$ has full rank it suffices to show that its null space is trivial. 

If $z\in\R^k$ is a vector in the null space of $\mA_1'$, then
\[
 \mA_1' z = \mP^\perp \mA_1 z = 0.
\]
If
\[
 \mA_1 z = \sum_{j=1}^k z_j a^{(j)},
\]
then
\[
 \mP^\perp\left(\sum_{j=1}^k z_j a^{(j)}\right) = 0,
\]
hence
 \[
  \sum_{j=1}^k z_j a^{(j)} \in {\mathcal N}(\mP^\perp) = {\mathcal R}(\mA_2)
 \]   
 and 
 \[
   \sum_{j=1}^k z_j a^{(j)}  \in  {\mathcal R}(\mA_1)\cap  {\mathcal R}(\mA_2) = \{0\}.
 \]
The linear independence of the columns of $\mA_1$ implies that $z = 0$, hence $\mA_1'$ is full rank.

The full rank of $\mA_1'$ implies that the projected equation in the spotlight,  (\ref{projected}) has the unique (least squares) solution
\[
 x_1' = \left(\mA_1'\right)^\dagger b' = \left((\mA_1')^\mT\mA_1'\right)^{-1}\left(\mA_1'\right)^\mT b'.
 \]
Let
\[
  x = \left[\begin{array}{c} x_1 \\ x_2 \end{array}\right] 
 \]
 be the unique solution of the original linear system (\ref{full}). Its orthogonal projection onto the orthogonal complement of the range of $\mA_
 2$ satisfies 
 \[
  \mA'_1 x_1 = b',
 \]
 and by the uniqueness of the solution it follows that  $x_1' = x_1$. \hfill $\Box$

\section{Spotlight projections for inverse problems}\label{sec:bayesian}

In this section, we consider the case in which the vector $b$ represents a noisy observation, and the matrix $\mA$ is ill-conditioned. Consequently, the inverse problem is ill-posed and therefore it requires  regularization, or alternatively, its stochastic extension can be considered in the Bayesian framework. The proposed approach is to consider the projected linear inverse problem
\[
 b' = \mA_1' x_1 + \varepsilon'
\]
as an ill-posed problem independent of the original one, introducing any standard regularization scheme, or, alternatively, a prior for the unknown $x_1$.
While the spotlight projection does not assume or require a Bayesian interpretation of the problem, it is natural to ask whether the methodology lends itself to a statistical interpretation. In the following subsection,  
we provide a Bayesian interpretation of spotlight projection method, under the assumption of Gaussian prior and Gaussian likelihood. We begin by outlining standard techniques to deal with nuisance parameters based on marginalization and modeling the clutter as part of the noise. For a general reference on the Bayesian solution of linear inverse problems, see, e.g., \cite{calvetti2023bayesian}.

\subsection{Marginalization and error modeling}

In the Bayesian extension of the inverse problem (\ref{original}), let $X$ denote a $\R^n$-valued Gaussian random variable with prior probability density $X\sim{\mathcal N}(0,\mC)$, where $n=n_1+n_2$, and $\mC\in\R^{n\times n}$ is a symmetric positive definite covariance matrix. We partition $X$ and the forward model 
\[
 X = \left[\begin{array}{c} X_1 \\ X_2\end{array}\right], \quad X_j\in\R^{n_j}, \quad \mA = \left[\begin{array}{cc} \mA_1 & \mA_2\end{array}\right] \in \R^{m\times n},
\]
and write the stochastic extension of the model as
\[
B = \mA X + E = \mA_1 X_1 + \mA_2 X_2 + E,
\]
where $E$ is an additive Gaussian noise vector independent of $X$. Without loss of generality, we may assume that the noise is whitened, i.e., $E\sim{\mathcal N}(0,\mI)$.  Furthermore, we partition the prior covariance matrix as
\[
 \mC = \left[\begin{array}{cc} \mC_{11} & \mC_{12} \\ \mC_{21} & \mC_{22}\end{array}\right], \quad \mC_{jk} \in\R^{n_j\times n_k}.
\] 

The estimation of the variable $X_1$ based on the observation $B$, may be 
done in two different ways:
\begin{enumerate}
\item {\em Lumping and conditioning:}  Write the forward model as
\[
 B = \mA_1 X_1 + \big(\mA_2 X_2 + E\big) = \mA_1 X_1 + \widehat E,
\]
interpreting the clutter term as part of the observation noise $\widehat E = \mA_2 X_2 + E$, and then calculate directly the posterior density $\pi_{X_1\mid B}(x_1\mid b)$ by conditioning. 
\item {\em Marginalization:} Calculate the posterior probability density of $(X_1, X_2)$ 
\[
 \pi_{X\mid B}(x\mid b) = \pi_{X_1,X_2\mid B}(x_1,x_2\mid b),
\]
and then marginalize it with respect to $X_2$,
\[
 \pi_{X_1\mid B}(x_1\mid b) = \int \pi_{X_1,X_2\mid B}(x_1,x_2\mid b) {\rm d} x_2. \label{eq:marginalX1}
\]  
\end{enumerate}

While it is intuitively clear that these two approaches should lead to the same posterior density, for completeness we state this as a theorem and give a short proof of it.
\begin{theorem}
Under the above hypotheses, the derivations of the posterior density $\pi_{X_1\mid b}(x_1\mid b)$ via conditioning and marginalization are equivalent. Moreover, the posterior is a Gaussian density,
\[
 X_1\mid B\sim{\mathcal N}(\mu_1,\mD_1),
\]
with mean and covariance 
\[
\mu_1 = \big(\mC_{11}\mA_1^\mT + \mC_{12}\mA_2^\mT \big)\big( \mA \mC \mA^\mT + \mI \big)^{-1} b, \quad 
\mD_1 = \mC_{11} - \big(\mC_{11}\mA_1^\mT + \mC_{12}\mA_2^\mT \big)\big( \mA \mC \mA^\mT + \mI \big)^{-1}\big( \mA_1 \mC_{11} + \mA_2\mC_{21} \big).
\] 
\end{theorem}

{\em Proof.}  
{\bf Conditioning approach.}   Combine $X_1$ and $B$ into the Gaussian random variable
\[
 Z = \left[\begin{array}{c} X_1 \\ B\end{array}\right].
\] 
%\[
% \bbE\big(Z Z^\mT\big) =\left[\begin{array}{cc} \bbE\big(X_1 X_1^\mT\big) & \bbE\big(X_1 B^\mT\big) \\   \bbE\big(B X_1^\mT\big) & \bbE\big(B  B^\mT\big) \end{array}\right]
%\] 
From the observations that   
\[
 \bbE\big(X_1 X_1^\mT\big) = \mC_{11}, \quad \bbE\big(X_1 B^\mT\big) = \mC_{11}\mA_1^\mT + \mC_{12}\mA_2^\mT = \bbE\big(B X_1^\mT\big)^\mT,
\]
and
\[
 \bbE\big(B B^\mT\big) =
  \bbE\big((\mA X + E) (\mA X + E)^\mT\big)
 = \mA \mC \mA^\mT + \mI,
\]

it follows that 
\[
  \bbE\big(Z Z^\mT\big) = \mD = \left[\begin{array}{cc} \mD_{11} & \mD_{12} \\ \mD_{21} & \mD_{22}\end{array}\right]\\
  = \left[\begin{array}{cc} \mC_{11} & \mC_{11}\mA_1^\mT + \mC_{12}\mA_2^\mT 
   \\ \mA_1 \mC_{11} + \mA_2\mC_{21} & \mA \mC \mA^\mT + \mI.\end{array}\right].
 \]  
Thus, the mean and the covariance matrix of $X_1$ conditional on $B=b$ can be expressed in terms of the Schur complement of $\mD_{22}$ as 
\begin{eqnarray*} 
\bbE(X_1\mid B = b) & = & \big(\mC_{11}\mA_1^\mT + \mC_{12}\mA_2^\mT \big)\big( \mA \mC \mA^\mT + \mI \big)^{-1} b, \\
 {\rm cov}(X_1\mid B = b) & = & \mC_{11} - \big(\mC_{11}\mA_1^\mT + \mC_{12}\mA_2^\mT \big)\big( \mA \mC \mA^\mT + \mI \big)^{-1}\big( \mA_1 \mC_{11} + \mA_2\mC_{21} \big),
\end{eqnarray*}
which agree with the expressions given in the theorem. 

{\bf Marginalization approach.}  The joint posterior density of the random variable $X = [ X_1; X_2]^\mT $ given $B=b$ is
\[
 \pi_{X_1,X_2\mid B}(x_1,x_2\mid b) = {\mathcal N}(x\mid \mu,\mD), \quad 
 \mu = \mC \mA^\mT\big(\mA\mC \mA^\mT + \mI\big)^{-1} b, \quad \mD =\mC - \mC \mA^\mT (\mA\mC\mA^\mT + \mI)^{-1}\mA\mC.
\] 
It can be shown, after completing the square and computing Schur complements, that the marginal density of $X_1$ is
\[
 \pi_{X_1\mid B}(x_1\mid b) = {\mathcal N}(x_1\mid \mu_1,\mD_{11}),
\]
where $\mu_1$ and $\mD_{11}$ refer to block partitionings of the posterior mean and covariance,
\[
 \mu= \left[\begin{array}{c} \mu_1 \\ \mu_2\end{array}\right], \quad \mD = \left[\begin{array}{cc} \mD_{11} & \mD_{12} \\ \mD_{21} & \mD_{22}\end{array}\right].
\] 
Recalling that 
\[
 \mC \mA^\mT = \left[\begin{array}{cc} \mC_{11} & \mC_{12} \\ \mC_{21} & \mC_{22}\end{array}\right]\left[\begin{array}{c} \mA_1^\mT \\ \mA_2^\mT\end{array}\right]
 = \left[\begin{array}{c} \mC_{11}\mA_1^\mT + \mC_{12}\mA_2^\mT \\ \mC_{21}\mA_1^\mT + \mC_{22}\mA_2^\mT\end{array}\right],
\]
we obtain
\[
 \mu_1 = \big( \mC_{11}\mA_1^\mT + \mC_{12}\mA_2^\mT\big)\big(\mA\mC \mA^\mT + \mI\big)^{-1} b,
\]
and
\[
 \mD_{11} = \mC_{11} - \big( \mC_{11}\mA_1^\mT + \mC_{12}\mA_2^\mT\big)\big(\mA\mC \mA^\mT + \mI\big)^{-1} \big( \mA_1\mC_{11} + \mA_2 \mC_{21}\big),
\]
thus we find that the density is the same as what we found in the previous part. \hfill$\Box$

\subsection{Projection versus marginalization}

The density of $X_1$ obtained through marginalization  differs from the posterior density of the projected problem (\ref{projected}) because the former depends on both $\mP b$ and $\mP^\perp b$. To shed some light on how the projection method relates to the to marginalization approach, we write
\begin{eqnarray*}
 \| b - \mA_1 x_1 - \mA_2 x_2\|^2 &=&  \| \mP^\perp(b - \mA_1 x_1 - \mA_2 x_2) +\mP(b - \mA_1 x_1 - \mA_2 x_2) \|^2  \\
  &=&   \|  b'  - \mA_1' x_1\|^2+ \| b'' - \mA_2 x_2 -  \mA_1'' x_1 \|^2 , 
 \end{eqnarray*}
 where
 \[
   b'' = \mP b, \quad \mA_1'' = \mP\mA_1.
 \] 
and substitute this decomposition into the whitened likelihood model $\pi_{B\mid X}(b\mid x) $ to obtain
 \begin{eqnarray}\label{likelihood factor}
 \pi_{B\mid X}(b\mid x) &\propto& {\rm exp}\left( -  \frac 12  \|  b'  - \mA_1' x_1\|^2\right) {\rm exp}\left( -\frac 12  \| b'' - \mA_2 x_2 -  \mA_1'' x_1 \|^2 \right)  \nonumber \\
  &\propto& \pi_{B'\mid X_1}( b'\mid x_1)\pi_{B''\mid X_1,X_2 }(b''\mid x_1,x_2). 
\end{eqnarray}
If, in addition, $X_1$ and $X_2$ are Gaussian and mutually independent, the joint prior factors as
\[
 \pi_X(x) = \pi_{X_1,X_2}(x_1,x_2) = \pi_{X_1}(x_1)\pi_{X_2}(x_2),
\]
and the posterior distribution can be written as 
\begin{eqnarray*}
\pi_{X_1,X_2\mid B} (x_1,x_2\mid b) 
&\propto& \pi_{B'\mid X_1}( b'\mid x_1)\pi_{B''\mid X_1,X_2}(b''\mid x_1,x_2) \pi_{X_1}(x_1)\pi_{X_2}(x_2)  \\
&=& \underbrace{\pi_{X_1}(x_1)\pi_{B'\mid X_1}(b'\mid x_1) }_{(a)} \; \underbrace{\pi_{X_2}(x_2)\pi_{B''\mid X_1,X_2}(b''\mid x_1,x_2)}_{(b)},
\end{eqnarray*} 
where (a) coincides with the posterior of the projected model. In general, since (b) depends on $x_1$, it cannot be ignored. The fact that marginalization with respect to $x_2$ will give a contribution to the posterior density of $X_1$, complicates the Bayesian interpretation of the projected problem in the general case.

The theorem below states that the posterior of the projected problem converges to the marginal density as the prior for $X_2$ becomes asymptotically uninformative. 
More specifically, assume that the prior for $X_2$ is of the form
\[
 \pi_{X_2}(x_2 ) =  \pi_{X_2}^\alpha(x_2 )\sim{\mathcal N}( 0, \alpha^{-2} \mG),
\]
where $\mG$ is a symmetric positive definite matrix, and $\alpha>0$ is a scaling parameter. Without loss of generality, possibly via a whitening transformation, we assume that $\mG = \mI$. As $\alpha\to 0$, the prior distribution of $X_2$ becomes increasingly wide, and the prior becomes asymptotically uninformative.  In the following theorem, we denote by $\pi^\alpha_{X_1,X_2\mid B}$ the conditional joint density when the prior for $X_2$ is chosen as a function of $\alpha$ in this manner. 

\begin{theorem}
Assume that $\mA_2\in\R^{m\times n_2}$ and ${\rm rank}(\mA_2) = n_2<m$. As the the prior for $X_2$ becomes uninformative, the marginal probability density of $X_1$, 
\[
 \pi^\alpha_{X_1\mid B'}(x_1\mid b) = \int \pi^\alpha_{X_1,X_2\mid B}(x_1,x_2\mid b) {\rm d} x_2,
\] 
converges to the posterior density for the projected problem (\ref{projected}), i.e., for some constant $C>0$,  
\[
 \lim_{\alpha\to 0} \alpha^{-n_2} \pi^\alpha_{X_1\mid B'}(x_1\mid b') = C \pi_{X_1}(x_1) \pi_{B' \mid X_1}(b'\mid x_1).
 \]
\end{theorem}

{\em Proof:}
We want to prove that in the limit $\alpha\to 0$, the integral
\[
 I_\alpha = \int \pi^\alpha_{X_2}(x_2)\pi_{B''\mid X_1,X_2}(b''\mid x_1,x_2) {\rm d} x_2,
\] 
is independent of $x_1$. To that end, we observe that
\begin{eqnarray*}
  \| b''- \mA_2 x_2 -  \mA_1'' x_1 \|^2 + \alpha^2 \|x_2\|^2 &=&
  x_2^\mT(\mA_2^\mT \mA_2 + \alpha^2 \mI)x_2 - 2x_2^\mT \mA_2^\mT(b'' - \mA_1'' x_1) + \|b''-\mA_1''x_1\|^2 \\
  &=& x_2^\mT \mH_\alpha x_2 - 2 x_2 \mA_2^\mT  y + \|y\|^2 \\
  &=& (x_2 - \mH_\alpha^{-1} \mA_2^\mT y)^\mT \mH_\alpha (x_2 - \mH_\alpha^{-1} \mA_2^\mT y) + \|y\|^2  - y^\mT \mA_2 \mH_\alpha^{-1}\mA_2^\mT y,
  \end{eqnarray*} 
  where
  \[
   \mH_\alpha = \mA_2^\mT \mA_2 + \alpha^2 \mI = \mH_0 +\alpha^2\mI, \quad y = b'' - \mA_1''x_1.
  \] 

 Since $\mA_2$ has full rank and $\mH_0$ is invertible,  we may use the resolvent identity
 \[
  \mH_\alpha^{-1} = \mH_0^{-1} - \alpha^2 \mH_0^{-1}\mH_\alpha^{-1},
 \]  
to write
 \[
   y^\mT \mA_2 \mH_\alpha^{-1}\mA_2^\mT y =  y^\mT \mA_2 \mH_0^{-1}\mA_2^\mT y -\alpha^2 y^\mT \mA_2 \mH_0^{-1}\mH_\alpha^{-1}\mA_2^\mT y.
 \]
Moreover, from $y\in{\mathcal R}(\mA_2)$ and the definition of $\mH_0$  we have  $\mA_2 \mH_0^{-1}\mA_2^\mT y = y$, and
 $ y^\mT \mA_2 \mH_0^{-1}\mA_2^\mT y = y^\mT y = \|y\|^2$, 
 implying that
 \[
  x_2^\mT \mH_\alpha x_2 - 2 x_2\mA_2^\mT  y + \|y\|^2 =  (x_2 - \mH_\alpha^{-1} \mA_2^\mT y)^\mT \mH_\alpha (x_2 - \mH_\alpha^{-1} \mA_2^\mT y) 
  +\alpha^2 y^\mT \mA_2 \mH_0^{-1}\mH_\alpha^{-1}\mA_2^\mT y.
 \]

 We can now evaluate the integral $I_\alpha$, 
 \begin{eqnarray*}
  I_\alpha &=& C'\alpha^{n_2}\int {\rm exp}\left(-\frac 12 \left(\| b'' - \mA_2 x_2 -  \mA_1'' x_1 \|^2 + \alpha^2 \|x_2\|^2\right)\right) {\rm d} x_2 \\
  %&=& C\alpha^{-n_2}\int {\rm exp}\left( -\frac 12(x_2 - \mH_\alpha^{-1} \mA_2^\mT y)^\mT \mH_\alpha (x_2 - \mH_\alpha^{-1} \mA_2^\mT y) \right) dx_2\, {\rm exp}\left(-\frac 12\alpha^2 y^\mT \mA_2 \mH_0^{-1}H_\alpha^{-1}\mA_2^\mT y\right) \\
  &=& \frac {(2 \pi)^{n_2/2}C' \alpha^{n_2}}{|{\rm det}(\mH_\alpha)|^{1/2}}{\rm exp}\left(-\frac 12\alpha^2 y^\mT \mA_2 \mH_0^{-1}\mH_\alpha^{-1}\mA_2^\mT y\right),
  \end{eqnarray*}
  where $C'$ is a constant, and hence
  \[
  \alpha^{-n_2} I_\alpha 
  \rightarrow \frac {(2 \pi)^{n_2/2}C'}{|{\rm det}(\mH_0)|^{1/2}} = C, \quad\mbox{as $\alpha\to 0$.}
  \]
  This completes the proof.
  \hfill$\Box$

For simplicity, in the theorem we assumed that the matrix $\mA_2\in\R^{m\times n_2}$ has full rank. If ${\rm rank}(\mA_2) = r<n_2$, then it follows from the lean SVD of $\mA_2$, 
 \[
\mA_2 = \sum_{j=1}^r \sigma_j u_j(v_j)^\mT = \mU_r\mD_r \mV_r, 
 \]
 that the clutter term is in an $r-$dimensional subspace of $\R^m$
 \[
  \mA_2 x_2 = \mA_2 x_2', \quad x_2' = \mU_r \mU_r^\mT x_2,
 \]
and the projectors are  
\[
 \mP = \mU_r\mU_r^\mT,\quad \mP^\perp = \mI - \mP.
\] 
Hence, we may replace $n_2$ by $r$.
In particular, this calculation implies that the nuisance parameter in the model is effectively $r$-dimensional.

\subsection{Partial projection, ill-determined rank and underdetermined problems}\label{sec:partial}

In the previous section, it was assumed that the dimensionality of the nuisance parameter does not exceed the dimension of the data space, and that the columns of the matrix $\mA_2$ are linearly independent, conditions which may not always be satisfied in actual applications. If the rank of the matrix $\mA_2$ is ill-determined, i.e., some of its singular values are close to the computational precision, determining its effective rank is an ill-posed problem and the projectors may be very sensitive to roundoff errors.  Furthermore, if the dimension of the nuisance parameter vector $n_2$ exceeds $m$, we may have ${\mathcal N}(\mA_2^\mT) = \{0\}$ implying that $\mP^\perp = 0$ and the right hand side of the projected problem does not carry any information about $x_1$.

In this section we illustrate how the projection method can be modified to address these cases. Let
 \[
  \mA_2 = \mU \mD \mV^\mT
 \]
 be a complete SVD of $\mA_2$.
 If $n_2<m$ and $\mA_2$ is full rank, as we assumed in the previous section,
 \[
  \lambda_1\geq \ldots \geq   \lambda_{n_2}>0, \quad 
  {\mathcal R}(\mA_2) = {\rm span}\{u_1,\ldots,u_{n_2}\},
 \] 
 where $u_j\in\R^m$ is the $j$th column of $\mU$, and hence $\mP = \mU_{n_2}\mU_{n_2}^\mT$, with $\mU_{n_2}= [ u_1, \ldots,u_{n_2}]$.
 
 Replacing the exact projector $\mP$
 in the projected equation by its rank $r$ approximation $\mP_r = \sum_{j=1}^r u_j u_j^\mT$ for  some $r\leq n_2$  yields 
 \[
  b_r' = \mP_r^\perp b =
  \mA_{1,r}' x_1 + \mP_r^\perp\mA_2 x_2 + \varepsilon_r',
 \]
 where
 \[
  \mA_{1,r}' = \mP_r^\perp\mA_1, \quad \varepsilon_r' = \mP_r^\perp\varepsilon, \quad  \mP_r^\perp = \mI - \mP_r = \sum_{j=r+1}^m u_j u_j^\mT.
 \]
 Observe that since
 \begin{equation}\label{PrA}
  \mP_r^\perp\mA_2 = \sum_{j=r+1}^{n_2}\lambda_j u_j v_j^\mT,
 \end{equation}
 the projected clutter term vanishes if $r = n_2$. To estimate the expected size of the projected clutter we introduce the symmetric positive definite covariance matrix $\mC_{22}$ of the random variable $X_2$, the stochastic extension of the unknown vector $x_2$. Hence,
 \begin{eqnarray*}
\bbE\big(\|\mP_r^\perp\mA_2 X_2\|^2\big) &=&\bbE\big({\rm Trace}(\mP_r^\perp\mA_2 X_2)(\mP_r^\perp\mA_2 X_2)^\mT\big) \\
&=&
{\rm Trace}\big(\mP_r^\perp\mA_2 \bbE(X_2 X_2^\mT) (\mP_r^\perp\mA_2)^\mT\big)  \\
&=& {\rm Trace}\big(\mP_r^\perp\mA_2 \mC_{22}(\mP_r^\perp\mA_2)^\mT\big) \\
%&=&{\rm Trace}\big((\mP_r^\perp\mA_2 \mC_{22}^{1/2})(\mP_2^\perp\mA_2\mC_{22}^{1/2})^\mT\big) \\
&=&\big\|\mP_r^\perp\mA_2 \mC_{22}^{1/2}\big\|_F^2 = \big\|\sum_{j=r+1}^{n_2} \lambda_j u_j \big(\mC_{22}^{1/2} v_j\big)^\mT \big\|_F^2.
 \end{eqnarray*}
Furthermore, it follows from the orthonormality of the vectors $u_j$  that %the Frobenius norm of the sum of the rank-1 matrices can be expressed as
 \[
 \big\|\sum_{j=r+1}^{n_2} \lambda_j u_j \big(\mC_{22}^{1/2} v_j\big)^\mT \big\|_F^2 =
 \sum_{j=r+1}^{n_2} \lambda_j^2 \|\mC_{22}^{1/2} v_j\|^2,
 \]
 hence
  \[
 \bbE\big(\|\mP_r^\perp\mA_2 X_2\|^2\big) \leq
 \|\mC_{22}\|  \sum_{j=r+1}^{n_2} \lambda_j^2.
 \]
 We proceed in a similar manner to estimate the expected size of the projected noise $\varepsilon'_r$. Letting $E$ be the  whitened stochastic extension of  $\varepsilon$, $E \sim{\mathcal N}(0, \mI_m)$, we have that 
 \[
  \bbE\big(\|\mP_r^\perp E\|^2\big)
  = {\rm Trace}\big(\mP_r^\perp \mI_m (\mP_r^\perp)^\mT\big)
  ={\rm Trace}\big(\sum_{j={r+1}}^m u_j u_j^\mT\big) = m-r.
 \]

 Selecting the rank of the approximate projector as the smallest integer $r$ for which the condition
 \begin{equation}\label{r ineq}
   \|\mC_{22}\|  \sum_{j=r+1}^{n_2} \lambda_j^2 <m-r,
 \end{equation}
 is satisfied ensures that the expected size of the residual clutter is smaller than the projected additive noise, thus the projected clutter term can be safely neglected in the projected problem.  To have a computationally useful criterion for the truncation, we introduce the ratios
 \begin{equation}\label{eq:Rcurve}
  R_r = \frac{\|\mC_{22}\|  \sum_{j=r+1}^{n_2} \lambda_j^2}{m-r}, \quad 0\leq r <m.
 \end{equation}
 Observe that for $r=0$,
 the ratio 
 \[
  {\rm cSNR} = \frac{\bbE\big(\|\mA_2 X_2\|^2\big)}{\bbE\big(\|E\|^2\big)}
 \]
 represents the clutter signal-to-noise ratio, and we have the estimate
 \[
  {\rm cSNR} \leq  \frac{\|\mC_{22}\|  \sum_{j=1}^{n_2} \lambda_j^2}{m} = R_0. 
 \]
When the upper bound $R_0$ on the right in the above inequality is less than one, we may conclude that the clutter is dominated by the additive noise, and no projection is necessary. If, on the other hand, $R_0>1$, the ratios $R_r$ serve as a good indicator for selecting $r$. Observe that for $n_2<m$, we have $\lim_{r\to n_2} R_r = 0$, and the truncation value $r$ can be chosen as the smallest $r$ for which
\begin{equation}\label{eq:Rcriterion}
   R_r   <1, 
\end{equation}
yielding an estimate for the the smallest $r$ for which the inequality 
(\ref{r ineq}) is satisfied. For $n_2\geq m$ there is no guarantee that the inequality $R_r<1$ holds for any $r$, however, if it does, a feasible truncation value is the smallest $r$ for which it holds.  More generally, the value $r$ that minimizes $R_r$ gives the optimal clutter reduction.

Finally, we point out that the SVD of the matrix $\mA_2$ may be unfeasible in practice if the matrix dimensions are large. This is typically the case when the spotlight inversion technique is applied to inverse problems with spatially distributed parameters in three dimensions, the voxel number easily arriving to millions and the only feasible way to
%not even form the matrices explicitly but instead can 
operations by $\mA_2$ and $\mA_2^{\mT}$ is to carry out computations without explicit formation of $\mA_2$.
In such cases, a natural remedy is to use randomized linear algebra and algorithms, see, e.g.,
\cite{halko2011finding,martinsson2020randomized,tropp2023randomized} for results and further references. In particular, for a given integer $k$, one first forms a Gaussian random matrix
\[
 \mOmega \in \R^{n_2\times k}, \quad \mOmega_{ij}\sim{\mathcal N}(0,1).
\]
Next, the matrix $\mA_2$ is applied to the columns of $\mOmega$, yielding the matrix $\mZ = \mA_2\mOmega$. It can be shown that with high probability, the range of the matrix $\mZ$ gives a good approximation of the range of $\mA_2$, and an orthonormal basis for the range of $\mZ$ can be found with a relatively lightweight lean QR algorithm, or by calculating the lean SVD of $\mZ$, discarding singular values below the computational precision. We refer to the cited articles for more detailed discussion.

\subsection{Effect of Spotlighting on the Signal-to-Noise Ratio}

If the ranges of $\mA_1$ and $\mA_2$ are mutually orthogonal, the space $\R^m$ can be expressed as an orthogonal direct sum of these subspaces, and the spotlight projection is an orthogonal projection that does not affect the mapping $\mA_2$, i.e., $\mA_1' = \mP^\perp \mA_1 = \mA_1$. This is the ideal case for the spotlight inversion. In general, however, we expect a loss of part of the signal due to the oblique projection. In this subsection we analyze how the spotlight projection affects the signal-to-noise ratio.

Consider the original problem of estimating $x_1$ from the full observation model
\[
 b = \mA_1 x_1 + \big(\mA_2 x_2 + \varepsilon\big) = \mA_1 x_1 + \widetilde \varepsilon,
\]
where $\widetilde \varepsilon$ comprises both the observation noise and the clutter. To estimate the signal-to-noise ratio (SNR) of this model, we assume that, a priori, $x_1$ and $x_2$ are realizations of uncorrelated Gaussian random variables $X_1\sim{\mathcal N}(0,\mD_1)$ and $X_2\sim{\mathcal N}(0, \mD_2)$, and that $\varepsilon$ is scaled white noise, $E\sim{\mathcal N}(0,\sigma^2 \mI_m)$.
We write the expected power of the signal as
\[
 {\mathbb E}\big(\|\mA_1 X_1\|^2\big) = {\rm Trace}\left({\mathbb E}\big(\mA_1 X_1 X_1^\mT \mA_1^\mT\big)\right)
= {\rm Trace}\left(\mA_1 \mD_1  \mA_1^\mT\right) =\|\mA_1 \mD_1^{1/2}\|_F^2,
\]
and similarly, assuming that $X_2$ is independent of the observation noise, we write the expected power of the noise as
\[
{\mathbb E}\left(\|\widetilde\varepsilon\|^2\right) = \|\mA_2 \mD_2\|^2_F + m\sigma^2,
\]
thus we can estimate the signal-to-noise ratio as
\[
 {\rm SNR} = \frac{ {\mathbb E}\big(\|\mA_1 X_1\|^2\big)}
 {{\mathbb E}\left(\|\widetilde\varepsilon\|^2\right)} = \frac{\|\mA_1\mD_1^{1/2}\|^2_F}
 {\|\mA_2 \mD_2^{1/2}\|^2_F + m\sigma^2}.
\]

For comparison, consider the signal-to-noise ratio of the projected problem.
Let $\theta$ be the angle  between the subspaces ${\mathcal R}(\mA_1)$ and  ${\mathcal R}(\mA_2)$, defined as
\[
 \cos \theta = \sup\left\{|u^\mT v|,\; u \in {\mathcal R(\mA_1),\, v\in {\mathcal R}(\mA_2)},\; \|u\| = \|v\| = 1\right\}.
\]
We assume that $\theta>0$ and define 
\[
 \cos\theta = \sqrt{1 - \delta^2}, \quad \delta>0.
\]
In particular, if $\mP:\R^m\to{\mathcal R}(\mA_2)$ is the orthogonal projection, then for every $u\in{\mathcal R}(\mA_1)$, we have
\[
\|\mP u\|^2 =|u^\mT \mP u| \leq
\|u\|\|\mP u\||\cos\angle(u,\mP u)|
\leq \|u\|\|\mP u\||\cos \theta,
\]
hence
\[
 \|\mP u\|  \leq \|u\| \sqrt{1-\delta^2}.
\]
Therefore, the orthogonal projection $\mP^\perp$ satisfies
\[
 \|\mP^\perp u\|^2 = \| u\|^2 - \|\mP u\|^2 \geq \delta^2 \|u\|^2.
\]
To estimate the SNR of the projected problem we write
\[
{\mathbb E}\left(\|\mA_1' X_1\|^2\right) 
=\|\mA_1'\mD_1^{1/2}\|_F^2.
\]
By denoting the columns of the matrix $\mA_1\mD_1^{1/2}$ by $\widetilde a^{(1)}, \ldots, \widetilde a^{(n_1)} \in {\mathcal R}(\mA_1)$, we obtain an estimate
\[
 \|\mA_1'\mD_1^{1/2}\|_F^2
  = \sum_{j=1}^{k} \|\mP^\perp \widetilde a^{(j)}\|^2\geq \delta^2 \sum_{j=1}^{k} \| \widetilde a^{(j)}\|^2 = 
  \delta^2 \|\mA_1\mD_1^{1/2}\|^2_F.
\]
On the other hand, the power of the projected noise becomes
\[
 {\mathcal E}\big(\|\mP^\perp E\|^2) = \sigma^2(m-r).
\]
We therefore see that after the spotlight projection, the SNR has an estimate
\[
 {\rm SNR}' = \frac{{\mathbb E}\left(\|\mA_1' X_1\|^2\right) }
 {{\mathcal E}\big(\|\mP^\perp E\|^2)}\geq \delta^2 \frac{\|\mA_1\mD_1^{1/2}\|^2}{\sigma^2(m-r)},
\]
allowing a direct comparison with the original SNR, as stated by the following theorem.

\begin{theorem}
 Assume that the angle $\theta$ between the subspaces ${\mathcal R}(\mA_1)$ and ${\mathcal R}(\mA_2)$  satisfies $\cos\theta = \sqrt{1-\delta^2}$ with $\delta>0$, that  $X_1$ and $X_2$ are mutually independent Gaussian random variables and the observation noise is scaled white noise with noise variance $\sigma^2$. Then the signal-to-noise ratio ${\rm SNR}'$ and the corresponding signal-to-noise-ratio ${\rm SNR}$ with the clutter term included in the noise
 satisfy
 \[
  {\rm SNR}'\geq \delta^2
  \frac{\|\mA_2\mD_2^{1/2}\|_F^2 +m\sigma^2}{(m-r)\sigma^2} {\rm SNR}.
 \]
\end{theorem}

\section{Computed examples}

In this section, we illustrate the viability of the spotlight approach by applying it to two different problems. The first example performs a  model reduction by spotlighting a portion of the computed solution of an elliptic PDE discretized by the finite element method, and the second example applies spotlight inversion in local tomography problem with noisy input data.

\subsection{Forward model of the electrical impedance tomography}

As the first application of the method with noiseless data ($\varepsilon = 0$), we consider the forward problem of the electrical impedance tomography. Let $\Omega\subset\R^2$ be a bounded set with a connected complement, modeling an object with conductivity distribution  $\sigma:\Omega \to \R$. We assume that $0<\sigma_m\leq \sigma(x)\leq  \sigma_M<\infty$ for some positive constants $\sigma_m$ and $\sigma_M$.  On the boundary $\partial\Omega$, $L$ contact electrodes are attached,  modeled as non-overlapping intervals $E_\ell\subset\partial\Omega$ of the boundary curve, $1\leq\ell\leq L$. 
We assume that electric currents $I_\ell$ are injected through the electrodes, satisfying the Kirchhoff condition
\begin{equation}\label{Kirchhoff}
 \sum_{\ell = 1}^L I_\ell = 0.
\end{equation} 
The induced static electric voltage potential $u:\Omega\to \R$ satisfies the elliptic equation
\begin{equation}\label{voltage eq}
\nabla\cdot\big(\sigma\nabla u\big) = 0 \mbox{ in $\Omega$,}
\end{equation}
with the boundary conditions
\begin{equation}\label{boundary 1}
 \int_{E_\ell} \sigma\frac{\partial u}{\partial n} dS = I_\ell, \quad  \sigma\frac{\partial u}{\partial n} \bigg|_{\partial\Omega\setminus \cup E_\ell} = 0.
\end{equation}
Since this Neumann type boundary condition is not sufficient to determine uniquely $u$ we complement it by assuming that each electrode has a characteristic contact impedance $z_\ell >0$. We impose the additional condition
\begin{equation}\label{boundary 2}
 \left(u + z_\ell\sigma\frac{\partial u}{\partial n}\right)\bigg|_{{E}_\ell} = U_\ell,
\end{equation}
where the constants $U_\ell$ are the electrode voltages that satisfy a ground condition, 
\begin{equation}\label{ground} 
 \sum_{\ell = 1}^L U_\ell = 0.
\end{equation}  

We consider the forward problem of solving for the pair $(u, U)$ when  $\sigma$, $z$ and $I$ are given.
In this setting, $U$, $I$ and $z$ are vectors in $\R^L$ with components $U_\ell$, $I_\ell$ and $z_\ell$ appearing in the model above.  
This forward model is referred to as the complete electrode model (CEM) \cite{somersalo1992existence}. For an extension of the model with certain computational and theoretical advantages, we refer to \cite{hyvonen2017smoothened,darde2022contact}.

It is well known \cite{somersalo1992existence} that the CEM has a unique (weak) solution $(u,U)\in H^1(\Omega)\times \R^L_0$, where $\R^L_0$ is the set of real $L$-vectors satisfying the condition (\ref{ground}), and moreover, the solution satisfies the weak form equation,
\begin{equation}\label{weak form}
{\mathscr B}_{\sigma,z}\big((u,U),(v,V)\big) = \sum_{\ell=1}^L I_\ell V_\ell \quad\mbox{for all $(v,V) \in H^1(\Omega)\times \R^L_0$,}
\end{equation}
where the coercive quadratic form is given by
\[
{\mathscr B}_{\sigma,z}\big((u,U),(v,V)\big) = \int_\Omega \sigma\nabla u\cdot \nabla v \,dx + \sum_{\ell =1}^L \frac 1{z_\ell}\int_{{E}_\ell}(u-U_\ell)(v-V_\ell) dS .
\]
The above expression gives a natural finite element formulation of the forward model. Let ${\mathscr T}_h = \big\{K_\nu\big\}_{\nu = 1}^{n_t}$ denote a triangular tessellation of $\Omega$, such that $\overline\Omega_h = \cup \overline K_\nu$ is a polygonal approximation of the domain, where $h>0$ is a symbolic mesh size parameter. In the following, we do not distinguish between the domains $\Omega$ and $\Omega_h$.
Further, we denote by  $\{\psi_j\}_{j=1}^{n_v}$ a corresponding piecewise linear Lagrange basis, where $n_v$ is the number of vertices $p_k$ in the mesh, $\psi_j(p_k) = \delta_{jk}$.  We write an approximation of the voltage potential $u$ in $\Omega$,
\[
 u = \sum_{j=1}^{n_v} u_j\psi_j,
\]
and let $\{{\mathcal{E}}_\ell\}_{\ell=1}^{L-1}$ denote a basis of $\R^L_0$, such that
\[
 U = \sum_{\ell=1}^{L-1}\alpha_\ell {\mathcal{E}}_\ell, \quad I = \sum_{\ell=1}^{L-1}\beta_\ell {\mathcal{E}}_\ell.
\] 
 By letting 
\[
 \big(u,U\big)  = \sum_{j=1}^{n_v} u_j \underbrace{\big(\psi_j,0\big)}_{\overline\psi_j} + \sum_{\ell=1}^{L-1} \alpha_\ell \underbrace{\big(0, {\mathcal{E}}_\ell\big)}_{\overline\psi_{n_v+\ell}}  = \sum_{j=1}^{n_v} u_j \overline\psi_j + \sum_{\ell = 1}^{L-1} \alpha_\ell \overline\psi_{n_v+\ell},
\] 
the weak form equation (\ref{weak form}) can be written component-wise as
\[
  \sum_{j=1}^{n_v} {\mathscr B}_{\sigma,z}\big(\overline\psi_j,\overline\psi_k \big)  u_j +  \sum_{\ell=1}^{L-1} {\mathscr B}_{\sigma,z}\big(\overline\psi_{n_v+\ell},\overline\psi_k \big) \alpha_\ell = \sum_{\ell = 1}^{L-1}\langle \overline\psi_{n_v+\ell},\overline\psi_k\rangle \beta_\ell, \quad 1\leq k\leq n_v+L-1,
 \]
or, in matrix form as
\begin{equation}\label{matrix}
 \mA x = \left[\begin{array}{cc} \mA^{11} &  \mA^{12} \\ \mA^{21} & \mA^{22}\end{array}\right]\left[\begin{array}{c} u\\  \alpha \end{array}\right] = 
  \left[\begin{array}{c} 0\\  \beta \end{array}\right] = b ,
\end{equation}
where
\begin{eqnarray*}
 \mA_{jk}^{11}    &=& {\mathscr B}_{\sigma,z}\big(\overline\psi_j,\overline\psi_k \big),\quad 1\leq j,k\leq n_v, \\
 \mA_{j\ell}^{12} &=& {\mathscr B}_{\sigma,z}\big(\overline\psi_j,\overline\psi_{n_v+\ell} \big) = K^{21}_{\ell j},\quad 1\leq j\leq n_v, \; 1\leq \ell \leq L-1,\\
 \mA_{\ell \ell'}^{22} &=& {\mathscr B}_{\sigma,z}\big(\overline\psi_{n_v+\ell},\overline\psi_{n_v+\ell'} \big),\quad 1\leq \ell,\ell'\leq L-1.
 \end{eqnarray*}  
To find the voltages, we need to solve the linear system (\ref{matrix}), which fits the template for spotlighting $u$. 

\begin{figure}
\centerline{\includegraphics[width=8cm]{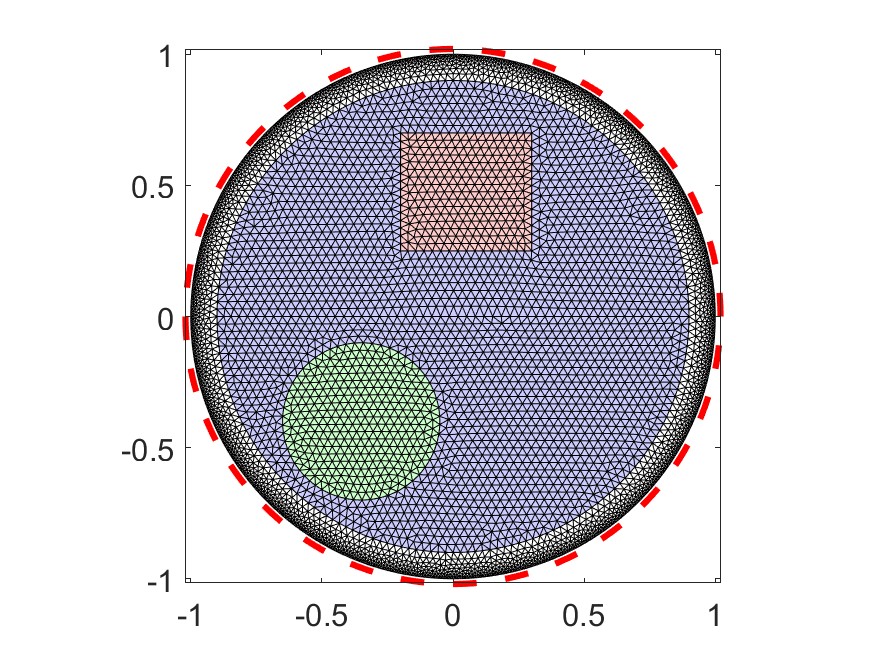}
\includegraphics[width=8cm]{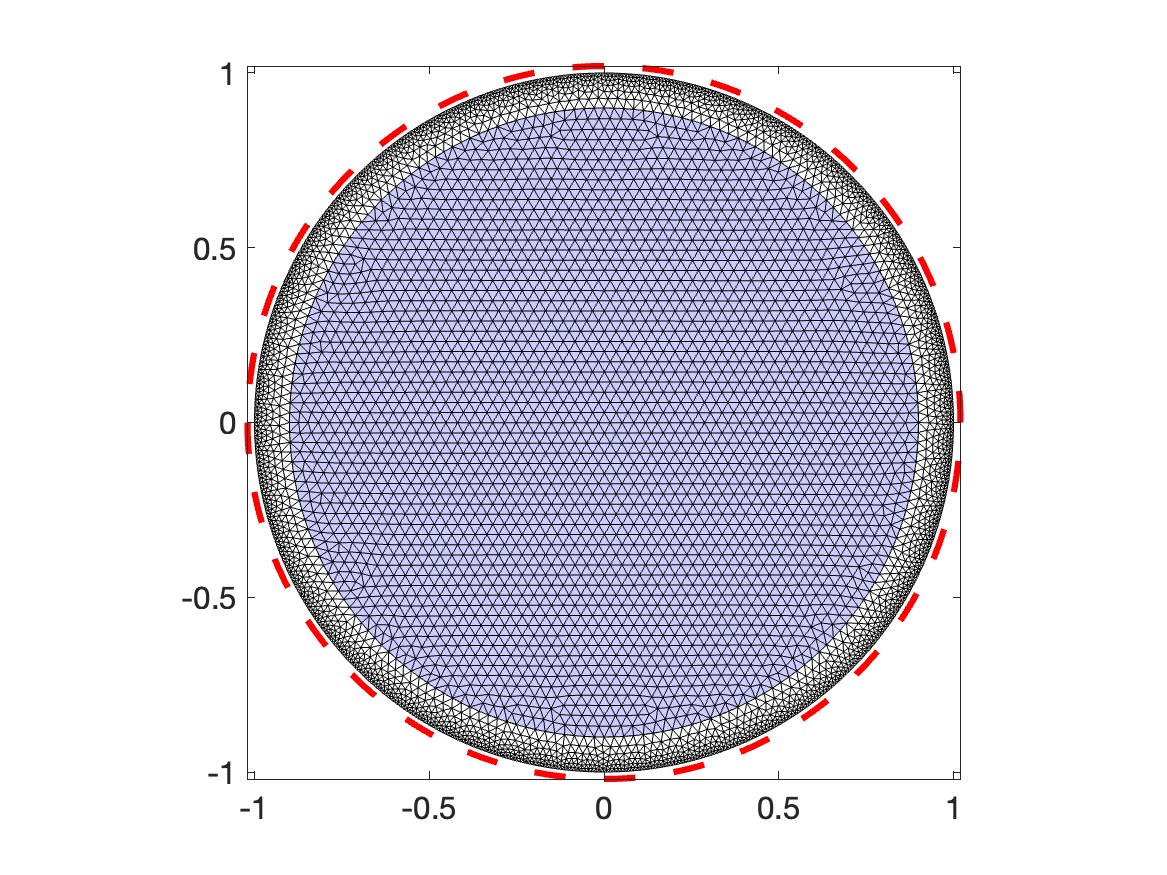}
}
\caption{Two FEM meshes with boundary refinement. The left one is used for data generation, the right one for solving the inverse problem. The 32 electrodes are indicated by the red line segments around the domain.}\label{fig:meshes}
\end{figure}

The boundary condition involving the electrodes creates a singularity in the voltage potential, the tangential derivative of the voltage along the electrodes increasing logarithmically as a function of the distance from the electrode edge. In order to properly account for this singularity, the FEM mesh needs to be very fine towards the boundary of the domain $\Omega$ in order to avoid boundary artifacts when solving the inverse problem. The boundary refinement typically increases the computational burden, even under the assumption that the conductivity near the boundary is known. Figure~\ref{fig:meshes} shows two FEM meshes on a unit disc with the boundary refinement performed in a collar region $\Omega_c = \{ x\in\Omega \mid |x|> 0.9\}$.   

To reduce the computational burden,  partition of the (ordered) index set ${\mathcal J} = \{1,2,\ldots,n_v\}$ into two subsets ${\mathcal J}_{\rm in}$ and ${\mathcal J}_{\rm out}$ 
\[
   {\mathcal J} = {\mathcal J}_{\rm in}\cup {\mathcal J}_{\rm out}, \quad {\mathcal J}_{\rm in}\cap {\mathcal J}_{\rm out} = \emptyset,
\]
setting
\[
 {\mathcal J}_{\rm in} = \big\{ k\in {\mathcal J}\mid |p_k|\leq 0.9\big\}, \quad     {\mathcal J}_{\rm out} = \big\{ k\in {\mathcal J}\mid |p_k| > 0.9\big\}.
\]  
Furthermore, we add to the index set ${\mathcal J}_{\rm in}$ with the indices pointing at the electrodes, 
\[
 {\mathcal J}_1 = {\mathcal J}_{\rm in}\cup\{n_v+1,\ldots,n_v+ L-1\}.
\] 
We define the matrices $\mA_1$ and $\mA_2$ as
\[
% \mA_1 = \mA(\,:\, ,{\mathcal J}_1), \quad  \mA_2 = \mA(\,:\, ,{\mathcal J}\setminus {\mathcal J}_1),
 \mA_1 = \mA(\,:\, ,{\mathcal J}_1), \quad  \mA_2 = \mA(\,:\, ,{\mathcal J}_{\rm out}),
\] 
The spotlight projection transforms the problem into a smaller one, where the unknowns are the voltage values in the interior domain $\Omega\setminus\Omega_c$,  as well as at the electrodes, discarding the many nodes in the collar domain.

\begin{figure}
\centerline{
\includegraphics[width=6cm]{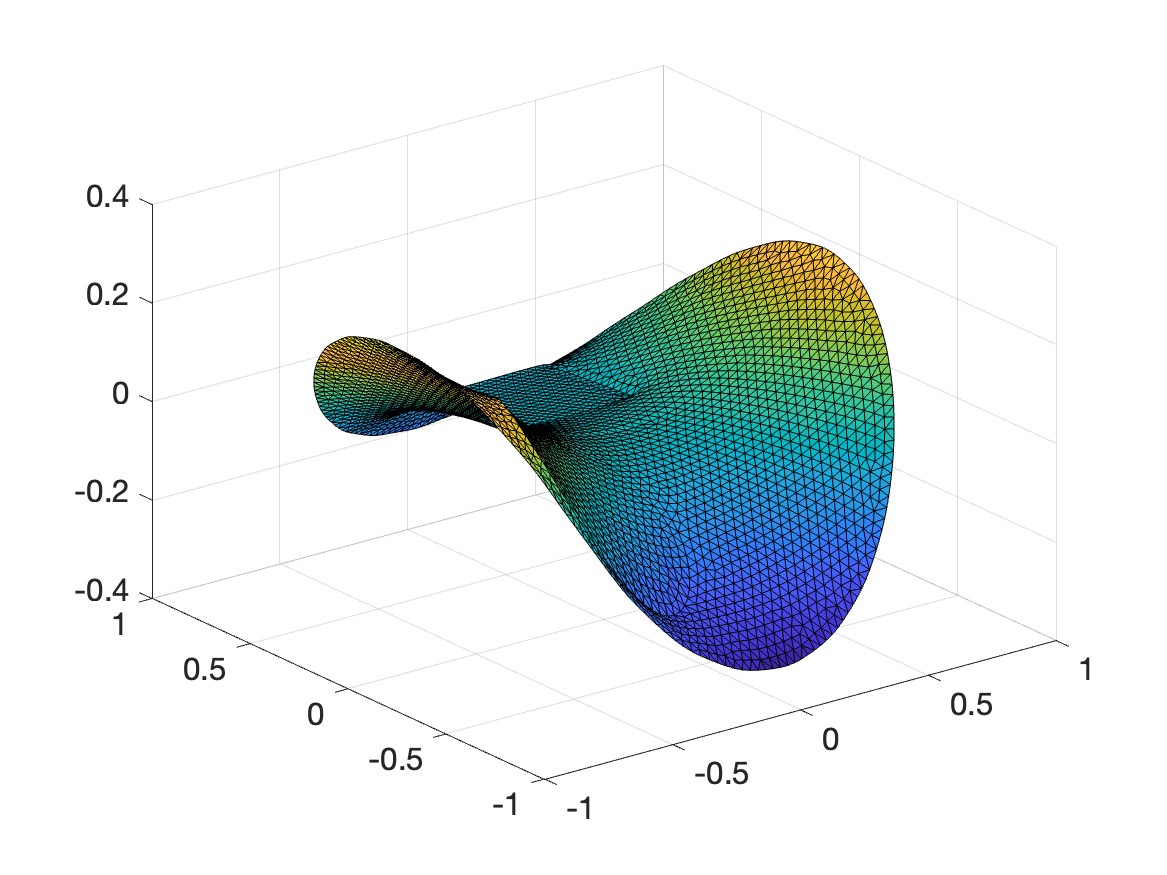}
\includegraphics[width=6cm]{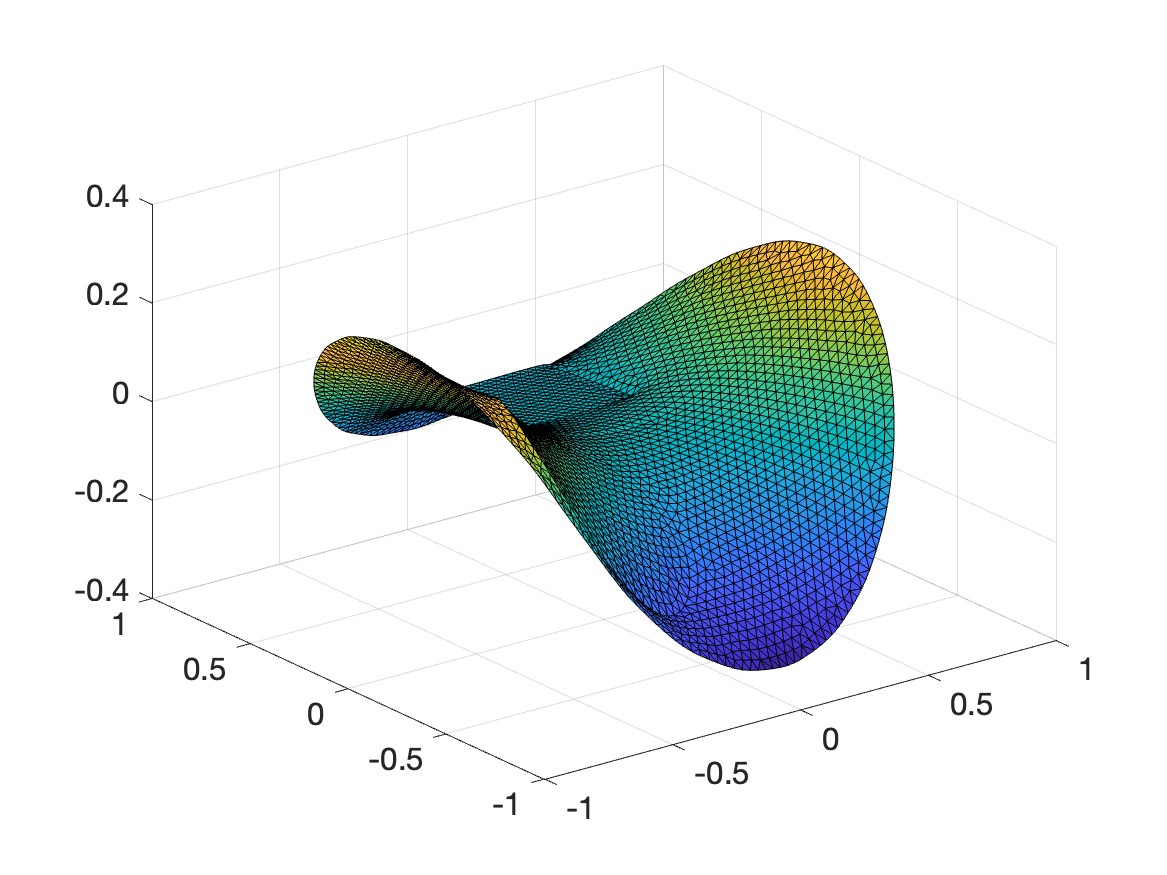}
\includegraphics[width=6cm]{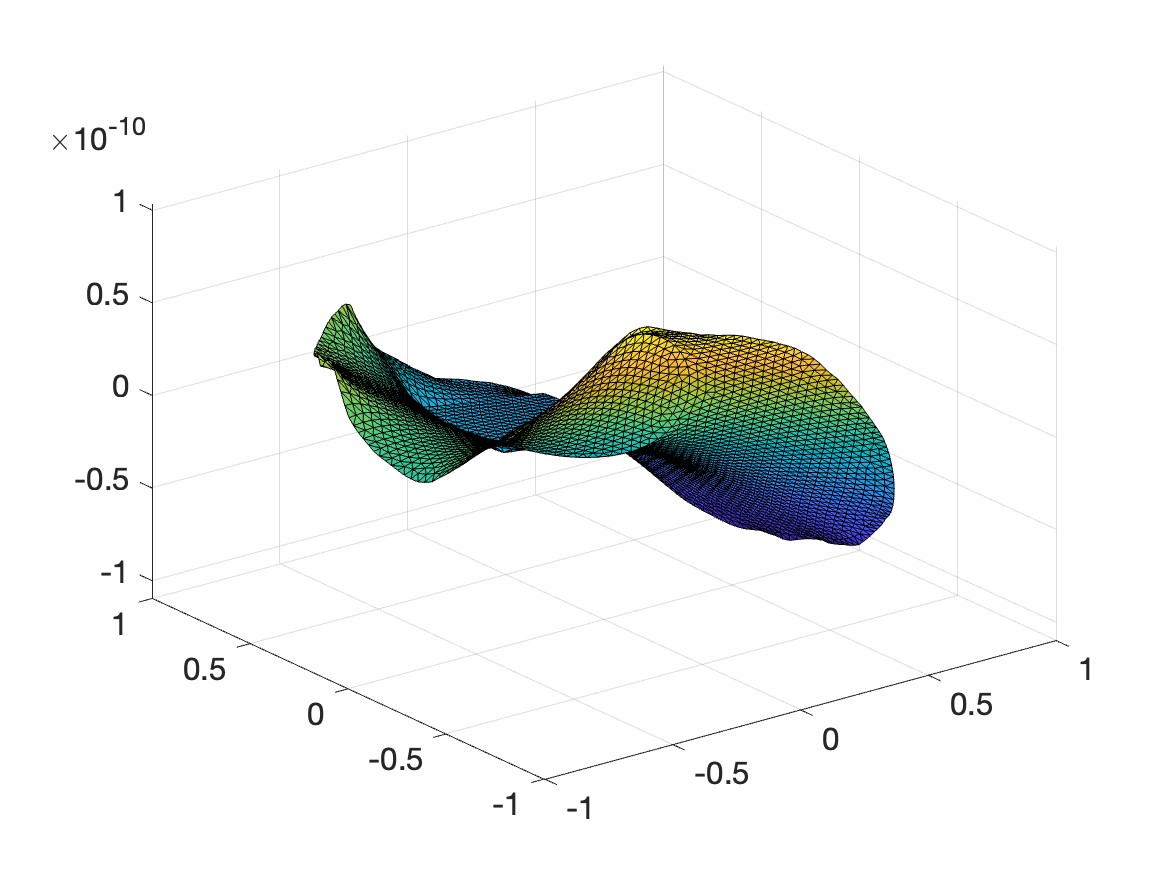}
}
\caption{Interior potential in $\Omega\setminus\Omega_c$ computed by using the full model (left) and the spotlight projections (middle). The difference, shown on the right, is negligible. }\label{fig:potentials}
\end{figure}

\begin{figure}
\centerline{
\includegraphics[width=10cm]{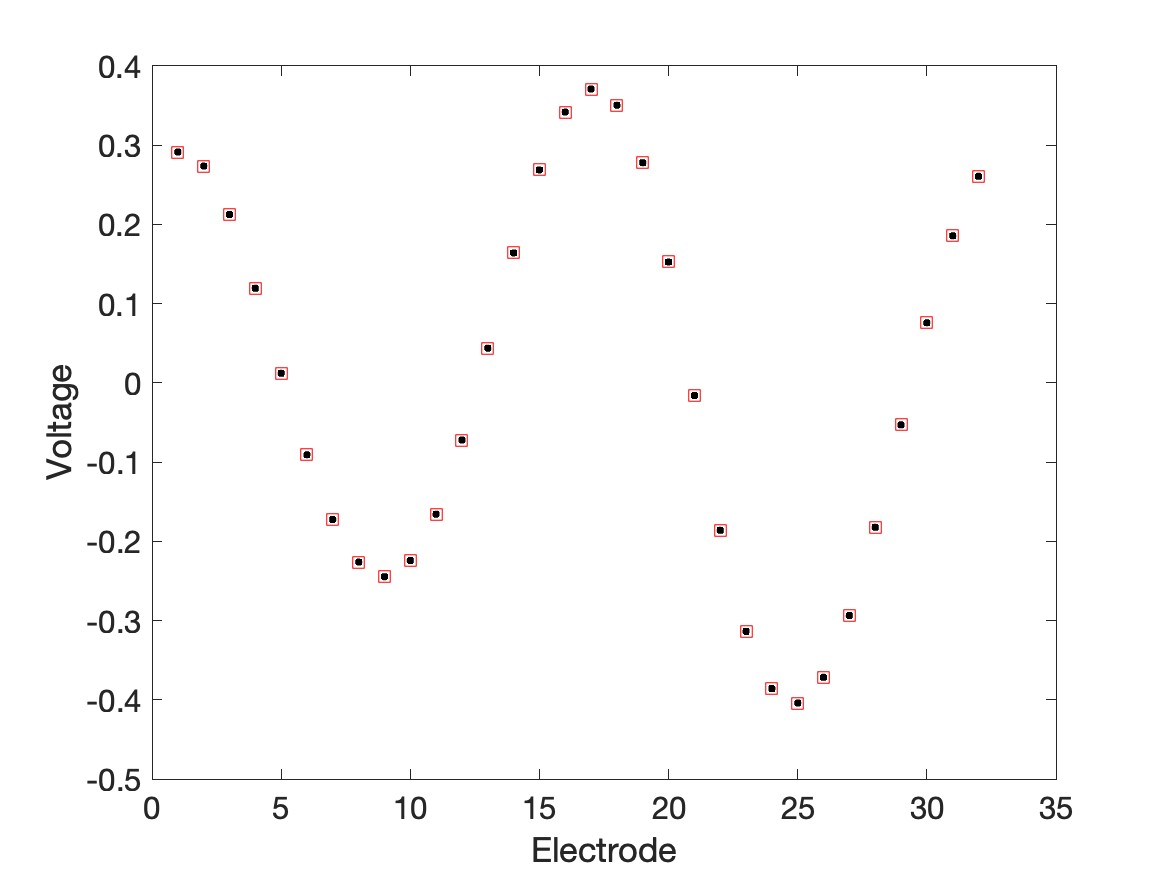}
}
\caption{Electrode voltages computed by the spotlight inversion (black dots), and the reference arising from the full model (red squares).}\label{fig:voltages}
\end{figure}

The voltage corresponding to a current injection $I_\ell = \cos 2 \frac{2\pi}{L} \ell$  obtained by using the spotlight inversion is shown in Figure~\ref{fig:potentials}. The full matrix is $\mA\in\R^{5842\times 5842}$, and the spotlight matrix is
$\mA_1 \in \R^{5842\times 2481}$, so we have roughly 42\%  of the unknowns left in the reduced  system. The comparison of the electrode voltages is  shown in Figure~\ref{fig:voltages}.

\subsection{Local tomography}

As a demonstration of the viability of the spotlight inversion approach for an ill-posed inverse problem, and to test the truncation criteria for the approximate projections, we apply the method to a computed two-dimensional local tomography problem with experimental fanbeam data. As a general reference to local tomography, see, e.g., the monograph \cite{ramm2020radon}. For a marginalization-based approach in the Bayesian inverse problems framework, we refer to
\cite{kolehmainen2011marginalization}.
 
The test case is based on an open access X-ray tomography data of a lotus plant root filled with different chemical elements \cite{bubbaOA}. The data set contains data for a 2D slice from fan beam scan of the target from 120 equally spaced angles over $360^\circ$ rotation, and forward matrix for the experiment using a $128 \times 128$ discretization of an image domain that covers completely the scanned specimen.    

Using this data, we formed a local tomography experiment by retrospective down-sampling of the data. We selected $40 \times 40$ pixel spotlight window (or region of interest, ROI) from the image domain and formed local tomography experiment by selecting the subset of X-rays that intersect the ROI, discarding all remaining data. This way, we ended having local tomography data $b \in \mathbb{R}^{14\,664}$ and forward model
\begin{equation}\label{eq:local_model}
b = \mA x + \varepsilon = \mA_1 x_1 + \mA_2 x_2 + \varepsilon
\end{equation}
with $x \in \mathbb{R}^{16\,384}$, $x_1 \in \mathbb{R}^{1\,600}$ and $x_2 \in \mathbb{R}^{14\,784}$. Observe that $n_2>m$, so a full spotlight projection is not feasible.

The level of measurement noise in the X-ray tomography experiment is not known. In the following, the noise $E$ is modeled as Gaussian scaled white noise, 
\begin{equation} \label{eq:noisemod}
    E \sim \mathcal{N}(0,\sigma^2 I) 
\end{equation}
with the noise standard deviation $\sigma$ estimated from a empty space patch (i.e., X-rays observed outside the attenuating target) of the original data, leading to an estimate for the standard deviation $\sigma = 0.0031$.  

\begin{figure}
    \centerline{
    \includegraphics[width=0.35\linewidth]{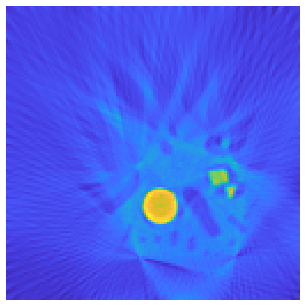} \quad
\includegraphics[width=0.35\linewidth]{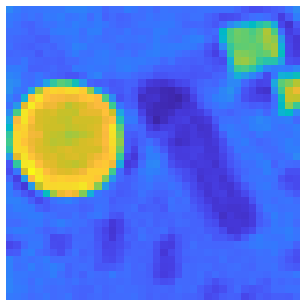}}
    \caption{Left: The MAP estimate from local tomography data using the full domain model $b=\mA x$. Right: Spotlight part from the MAP estimate on the left. We denote the spotlight part as $x^{\ast}_1$ and use it as ground truth reference for the estimates computed using the spotlight model $\mA_1 x_1$.}
    \label{fig:fullMAP}
\end{figure}

All computed reconstructions are Maximum A Posteriori (MAP) estimates, utilizing the noise model \eqref{eq:noisemod} and prior model 
\begin{equation}\label{eq:prior}
X \sim \mathcal{N}(0,\zeta^2 I).
\end{equation}
For the experiments, the prior standard deviation was set to $\zeta = 2.5 \cdot 10^{-4}$, and the feasibility of the selection was verified by visual assessment of the full domain local tomography reconstruction using the model $b= \mA x$. The MAP estimate using the full domain model is shown in Figure \ref{fig:fullMAP}. The spotlight part of the full domain model solution, shown on the right,
is the restriction of the solution to the ROI, and it serves as the ground truth for the estimates computed based on the spotlight model $\mA_1 x_1$. We denote the spotlight part of the MAP estimate with the model $\mA x$, or ground truth, as $x^{\ast}_1$. To assess the performance of the spotlight inversion, we define the relative error of any reconstruction $x_1$ by the formula
\begin{equation} \label{eq:image_error}
    d(x_1,x^{\ast}_1) = \frac{\| x_1 - x^{\ast}_1\|}{\| x^{\ast}_1\|}
\end{equation}

\begin{figure}
    \centering
    \includegraphics[width=0.7\linewidth]{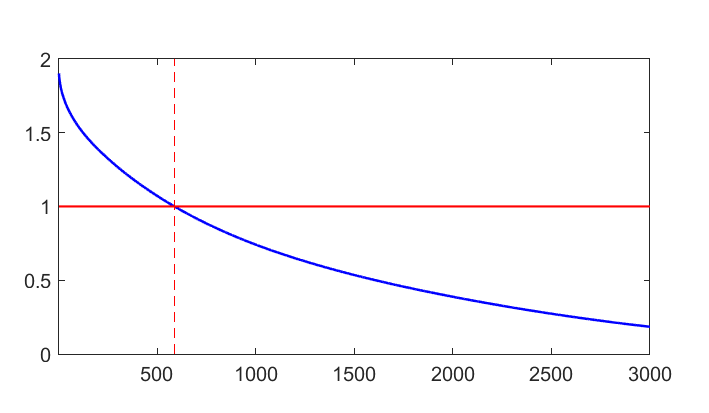}
    \caption{Clutter-to-noise ratio $R_r$ for the local tomography experiment. Horizontal axis is the number $r$ of svd vectors (i.e., number of columns in $U_r$). The point where $R_r =1$ is at $r=587$, denoted by the red dashed line.}
    \label{fig:Rcurve}
\end{figure}

To test the projection model, we computed SVD of $\mA_2$ and computed the clutter-to-noise ratio curve $R_r$ in equation \eqref{eq:Rcurve}, which is shown in Figure \ref{fig:Rcurve}. The clutter-to-noise ratio criterion
\eqref{eq:Rcriterion} yielded a truncation parameter $r=587$ for the approximate projection $
\mP_r^\perp = \mI - \mU_r \mU_r^{\mT}.
$

\begin{figure}
    \centerline{
    \includegraphics[width=0.35\linewidth]{Xmap_fullA_spotlight.png} \quad
    \includegraphics[width=0.35\linewidth]{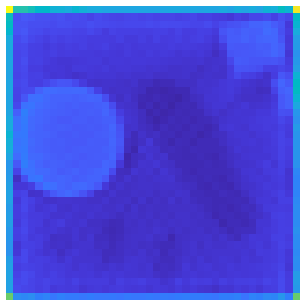}
    }
\vskip0.4cm
    
\centerline{
    \includegraphics[width=0.35\linewidth]{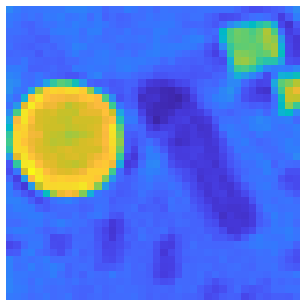}\quad
    \includegraphics[width=0.35\linewidth]{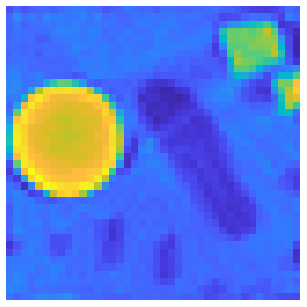}
    }
    
   \caption{Top left: Reference solution $x^{\ast}_1$ (detail from the MAP estimate with the full model $\mA x$). Top right: MAP estimate using the na\"{\i}ve model $b = \mA_1 x_1$ that ignores the clutter (relative error: $d(x_1,x^{\ast}_1)=1.51$). Bottom left: MAP estimate from marginal density $\pi_{X_1\mid B}(x_1 \mid b)$ (relative error: $d(x_1,x^{\ast}_1) \approx 10^{-15}$). Bottom right: MAP estimate using the projected spotlight model $\mP_r^\perp b = \mP_r^\perp \mA_1 x_1$ with $r=587$ (relative error: $d(x_1,x^{\ast}_1)=0.0596$). }
    \label{fig:result_R}
\end{figure}

Figure \ref{fig:result_R} shows the MAP estimates using the different models. For convenience, top left shows the reference image $x^{\ast}_1$ explained above. Top right shows the MAP estimate using the na\"{\i}ve reduced model $b= \mA_1 x_1$, assuming that contribution of the attenuation density distribution outside the spotlight domain is negligible, i.e., approximating $\mA_2 x_2 = 0$. As can be seen, the estimate exhibits large image artifacts, especially near the boundaries of the spotlight domain, and  it lacks in dynamical range, with the relative error compared to the reference $x^{\ast}_1$ (top left) around 150\%. The bottom left shows the MAP estimate from the marginal density $\pi_{X_1 \mid B} (x_1 \mid b)$, the mean $\mu_1$ of the marginal density of Theorem 1.1, which is, as expected, close to the reference $x^{\ast}_1$ from the MAP estimate with the full local tomography model $\mA x$. The bottom right shows the MAP estimate with the projected spotlight model $\mP_r^\perp b = \mP_r^\perp \mA_1 x_1$ with $r=587$, the truncation parameter value obtained by the cSNR analysis. The estimate is qualitatively very similar to the reference, with the relative error being around 6\%.

\begin{figure}
    \centering
    \includegraphics[width=0.7\linewidth]{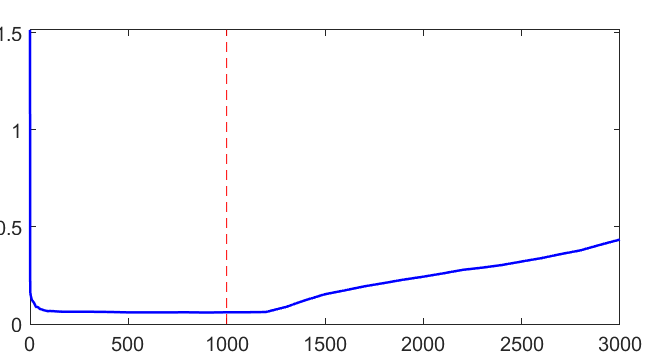}
    \caption{Estimate error $d(x_1,x_{\ast,1})$ with respect the truncation radius $r$ in the projection model. The minimum error is at $r=1000$ with $d(x_1,x_{\ast,1})=0.0595$}
    \label{fig:L2error_wrt_r}
\end{figure}

\begin{figure}
    \centerline{
     \includegraphics[width=0.35\linewidth]{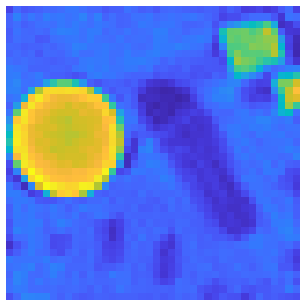}
   \quad \includegraphics[width=0.35\linewidth]{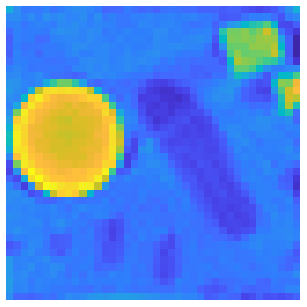}
   }
    \caption{Left: MAP estimate with the projected model $P_r^\perp b = P_r^\perp A_1 x_1$ with $r=1000$ corresponding to the minimum point in Figure \ref{fig:L2error_wrt_r} ($d(x_1,x_{\ast,1})=0.0595$). Right: MAP estimate with $r=100$ 
    (relative error: $d(x_1,x_{\ast,1})=0.0667$). }
    \label{fig:result_r100_and_rmin}
\end{figure}

To study the effect of the SVD truncation parameter $r$, we computed MAP estimates for a range of number of SVD  components (i.e. number of columns in $\mU_r$) included in the approximate projector. The estimation errors $d(x_1,x^{\ast}_1)$ with respect to truncation parameter $r$ are shown in Figure~\ref{fig:L2error_wrt_r}. 
We observe that the estimation error demonstrates a typical semi-convergence behavior, which can be understood through an SNR analysis: For small $r$, the clutter term dominates, while for high $r$, the projector $\mP_r^\perp$ starts to erode the information contents of the informative term $\mP_r^\perp \mA_1 x_1$, leading eventually to poor reconstruction.
The smallest estimation error occurs at $r=1\,000$ with $d(x_1,x^{\ast}_1) = 0.0595$, suggesting that the selection of the truncation parameter by using the clutter signal-to-noise ratio led to a highly feasible selection, as the difference in the relative error with $r=587$ was only 0.01\% compared to the minimum at $r=1\,000$. We observe from Figure \ref{fig:L2error_wrt_r} that the estimation error curve plateaus fast with respect to the number of SVD components, suggesting that already a relatively low number of singular vectors are enough for effective projection.
For comparison, the right image in Figure \ref{fig:result_r100_and_rmin} shows the MAP estimate with $r=100$, with the relative error being only slightly larger than with the estimate with the selection of $r$ based on the clutter signal-to-noise ratio (bottom right in Figure \ref{fig:result_R}) or in the estimate
with the smallest error compared to reference $x^{\ast}_1$ with
$r=1\,000$.  

\section{Conclusions and outlook}\label{sec:outlook}

This article proposes a purely linear algebraic approach for reducing a linear system with clutter, a linear contribution from a nuisance parameter, which, in spite of its limitations, is a commonly encountered type of clutter in practical problems.
The proposed spotlighting approach is demonstrated with two different numerical examples; the first example is related to solving the selected components from a finite element model of an elliptic PDE and the second to estimating the selected components in a local tomography problem.
The results with the computed examples indicate that this approach allows computationally efficient treatment of clutter, in particular in high resolution problems in three and four dimensions where possibly time dependency enters, and the use of full model or marginalization can be computationally expensive. 
In practice, and in light of some preliminary tests, we expect that for ill-posed inverse problems, one would need only few hundreds of left singular vectors of $\mA_2$ to reduce significantly the clutter, whereby the methods of randomized linear algebra will be of utmost utility.

The proposed projection method has direct connections to earlier works on nuisance parameter estimation. In fact, the formula (\ref{likelihood factor}) is a version of the likelihood factorization methods \cite{basu2010elimination}, taking advantage of the known fact that independency and orthogonality in the Gaussian framework coincide. This connection may provide new ideas of developing the more general nuisance parameter suppression methods towards an algorithmic direction.

The method has extensions to non-linear inverse problems, e.g., via applying spotlight inversion to equations obtained via sequential linearizations. Such an approach is particularly compelling if the range of the derivative with respect to the nuisance parameters is independent of the point at which it is evaluated; in fact, this study was inspired by such sensitivity analysis in  \cite{jaaskelainen2024projection} on electrical impedance tomography.

\section*{Acknowledgements}

The work of DC was partly supported by the NSF grants DMS 1951446 and  DMS 2513481, and that of  ES by the NSF grants DMS 2204618 and DMS 2513481. The work of NH and VK was supported by
Research Council of Finland (RCF) Flagship of Advanced Mathematics for Sensing, Imaging and Modeling grant (nos.~358944 and 359181), the RCF Centre of Excellence in Inverse Modelling and Imaging grant (nos.~353084 and 353081) and RCF grants (nos.~359433 and 359434).

\bibliographystyle{siam} 
\bibliography{Biblio}

\end{document}